\documentclass[11pt, final]{amsart}
\usepackage{amssymb}
\usepackage[mathscr]{eucal}
\usepackage{showkeys}

\theoremstyle{plain}
\newtheorem{theorem}{Theorem}

\newtheorem{lemma}[theorem]{Lemma}
\newtheorem{corollary}[theorem]{Corollary}

\theoremstyle{remark}

\def\A{\mathcal A}
\def\B{\mathcal B}
\def\C{\mathbb C}
\def\clrng{\overline{\operatorname{rng}}\,}
\def\la{\langle}
\def\N{\mathbb N}
\def\ot{\otimes}
\def\R{\mathbb R}
\renewcommand\Re{\operatorname{Re}}
\def\ra{\rangle}

\begin{document}
\title[]{Jordan maps on standard operator algebras}
\author{LAJOS MOLN\'AR}
\address{Institute of Mathematics and Informatics\\
         University of Debrecen\\
         4010 Debrecen, P.O.Box 12, Hungary}
\email{molnarl@math.klte.hu}
\thanks{  This research was supported by the
          Hungarian National Foundation for Scientific Research
          (OTKA), Grant No. T030082, T031995, and by
          the Ministry of Education, Hungary, Reg.
          No. FKFP 0349/2000}
\date{\today}
\begin{abstract}
Jordan isomorphisms of rings are defined by two
equations. The first one is the equation of additivity while
the second one concerns multiplicativity with respect to
the so-called Jordan product. In this paper we present results
showing that on standard operator algebras
over spaces with dimension at least 2,
the bijective solutions of that second equation are automatically
additive.
\end{abstract}
\maketitle

\section{Introduction and statement of the results}

It is an interesting problem to study the interrelation between the
multiplicative and the additive structures of a ring. The first quite
surprising result on how the multiplicative structure of a ring
determines
its additive structure is due to Martindale \cite{Martindale}. In
\cite[Corollary]{Martindale} he proved that every
bijective multiplicative
map from a prime ring containing a nontrivial idempotent onto
an arbitrary ring is necessarily additive and, hence, it is a ring
isomorphism.
This result has been utilized by \v Semrl in \cite{Semrl} to describe
the form of the semigroup isomorphisms of standard operator algebras on
Banach spaces.

Beyond ring homomorphisms there is another very important class of
transformations between rings. These are the Jordan
homomorphisms. The Jordan structure of associative rings has been
studied by many people in ring theory.
Moreover, Jordan operator
algebras have serious applications in the mathematical foundations of
quantum mechanics.
If $\mathcal R, \mathcal R'$ are rings and $\phi:\mathcal R\to
\mathcal R'$ is a transformation, then $\phi$ is called a Jordan
homomorphism if it is additive and satisfies
\begin{equation}\label{E:mjor14}
\phi(A^2)=\phi(A)^2 \qquad (A\in \mathcal R).
\end{equation}
If $\mathcal R'$ is 2-torsion free, then, under the assumption of
additivity, \eqref{E:mjor14} is equivalent to
\begin{equation}\label{E:mjor15}
\phi(AB+BA)=\phi(A)\phi(B)+\phi(B)\phi(A)
\qquad (A,B\in \mathcal R).
\end{equation}
Clearly, every ring homomorphism is a Jordan homomorphism and the same
is true for ring antihomomorphisms (the transformation $\phi:
\mathcal R \to \mathcal
R'$ is called a ring antihomomorphism if $\phi$ is additive and
satisfies $\phi(AB)=\phi(B)\phi(A)$ for all $A,B \in \mathcal R$).
In algebras, it seems more frequent that
instead of \eqref{E:mjor15},
one considers the equation
\begin{equation}\label{E:mjor16}
\phi((1/2)(AB+BA))=(1/2)(\phi(A)\phi(B)+\phi(B)\phi(A))
\qquad (A,B\in \mathcal R).
\end{equation}
Under the assumption of additivity these
two equations are obviously equivalent.

The additivity of bijective maps $\phi$ between von Neumann algebras
without commutative direct summand satisfying the equation
\eqref{E:mjor16} and
\begin{equation*}
\phi(A^*)=\phi(A)^*
\end{equation*}
was studied in \cite{Hakeda} and \cite{HakSai}.
The aim of this paper is to investigate similar problems on
standard operator algebras. From our present point of view the main
difference
between standard operator algebras and von Neumann algebras is
that standard operator algebras need not have unit and they are not
necessarily closed in any operator topology.

In what follows we shall study the equations \eqref{E:mjor15} and
\eqref{E:mjor16}. Both equations play important role; the first one is
because of ring theory and the second one is because of the applications of
Jordan operator algebras in mathematical physics.
Our main results describe the form of the bijective solutions of the
considered equations. It will turn out that all such solutions are
automatically additive.
We refer to our recent papers \cite{Molnar2}, \cite{Molnar}
for some other results of similar spirit.

We now summarize the results of the paper.
In what follows we consider all linear spaces over the complex field.
If $X$ is a Banach space, then we denote by $B(X)$ and $F(X)$ the
algebra of all bounded linear operators and the ideal of all bounded
linear finite rank operators on $X$, respectively. A subalgebra of
$B(X)$ is called a standard operator algebra if it contains $F(X)$.
The dual space of $X$ is denoted by $X'$ and $A'$ stands for the Banach
space adjoint of the operator $A\in B(X)$.

Our first result describes the form of the bijective solutions of
\eqref{E:mjor16} on standard operator algebras.

\begin{theorem}\label{T:mjor1}
Let $X, Y$ be Banach spaces, $\dim X>1$, and let $\A \subset B(X)$, $\B
\subset B(Y)$ be standard operator algebras.
Suppose that $\phi:\A \to \B$
is a bijective transformation satisfying
\begin{equation}\label{E:mjor1}
\phi((1/2)(AB+BA))=(1/2)(\phi(A)\phi(B)+\phi(B)\phi(A))
\end{equation}
for every $A,B \in \A$.

If $X$ is infinite dimensional, then we have the following
possibilities:
\begin{itemize}
\item[(i)]
there exists an invertible bounded linear operator $T: X\to Y$
such that
\[
\phi(A)=TAT^{-1} \qquad (A\in \A);
\]
\item[(ii)]
there exists an invertible bounded conjugate-linear operator $T: X\to
Y$ such that
\[
\phi(A)=TAT^{-1} \qquad (A\in \A);
\]
\item[(iii)]
there exists an invertible bounded linear operator $T: X'\to
Y$ such that
\[
\phi(A)=TA'T^{-1} \qquad (A\in \A);
\]
\item[(iv)]
there exists an invertible bounded conjugate-linear operator $T: X'\to
Y$ such that
\[
\phi(A)=TA'T^{-1} \qquad (A\in \A).
\]
\end{itemize}
If $X$ is finite dimensional, then we have $\dim X=\dim Y$.
So, our transformation $\phi$ can be supposed to act on the
matrix algebra $M_n(\C)$. In this case
we have the following possibilities:
\begin{itemize}
\item[(v)]
there exist a ring automorphism $h$ of $\C$ and an invertible matrix
$T\in M_n(\C)$ such that
\[
\phi(A)=Th(A)T^{-1} \qquad (A\in M_n(\C));
\]
\item[(vi)]
there exist a ring automorphism $h$ of $\C$ and an invertible matrix
$T\in M_n(\C)$ such that
\[
\phi(A)=Th(A)^tT^{-1} \qquad (A\in M_n(\C)).
\]
\end{itemize}
Here, ${}^t$ stands for the transpose and $h(A)$ denotes the matrix
obtained from $A$ by applying $h$ on every entry of it.
\end{theorem}

From the theorem above we easily have the following corollary.
If $H$ is a Hilbert space and $A\in B(H)$, then $A^*$ denotes the
Hilbert space adjoint of $A$.

\begin{corollary}\label{C:mjor3}
Let $H,K$ be Hilbert spaces, $\dim H>1$, and let $\A \subset B(H)$, $\B
\subset B(K)$ be
standard operator algebras which are closed under taking adjoints.
Suppose that $\phi:\A \to \B$
is a bijective transformation satisfying
\begin{equation}\label{E:mjor19}
\begin{aligned}
\phi(A^*) &=\phi(A)^*\\
\phi((1/2)(AB+BA)) &=(1/2)(\phi(A)\phi(B)+\phi(B)\phi(A))
\end{aligned}
\end{equation}
for every $A,B \in \A$.

Then we have the following possibilities:
\begin{itemize}
\item[(i)]
there exists a unitary operator $U: H\to K$
such that
\[
\phi(A)=UAU^* \qquad (A\in \A);
\]
\item[(ii)]
there exists an antiunitary operator $U: H\to K$ such that
\[
\phi(A)=UAU^* \qquad (A\in \A);
\]
\item[(iii)]
there exists a unitary operator $U: H\to
K$ such that
\[
\phi(A)=UA^*U^* \qquad (A\in \A);
\]
\item[(iv)]
there exists an antiunitary operator $U: H\to K$ such that
\[
\phi(A)=UA^*U^* \qquad (A\in \A).
\]
\end{itemize}
\end{corollary}

Unfortunately, we do not have a result concerning the equation
\eqref{E:mjor15} in the Banach space setting. However, we have
the following result describing the self-adjoint solutions of that
equation on standard operator algebras over Hilbert spaces.
The restriction to self-adjoint solutions is very
natural in operator theory where they
usually consider transformations on operator algebras which preserve
adjoints.

\begin{theorem}\label{T:mjor2}
Let $H,K$ be Hilbert spaces, $\dim H>1$, and let $\A \subset B(H)$, $\B
\subset B(K)$
be standard operator algebras which are closed under taking adjoints.
Let $\phi:\A \to \B$ be a bijective transformation satisfying
\begin{equation*}
\begin{aligned}
\phi(A^*) &=\phi(A)^* \\
\phi(AB+BA) &=\phi(A)\phi(B)+\phi(B)\phi(A)
\end{aligned}
\end{equation*}
for every $A,B \in \A$.
Then we have the same possibilities for $\phi$ as in
Corollary~\ref{C:mjor3}.
\end{theorem}

Although we do not have a result on the equation \eqref{E:mjor15} in the
Banach space setting, we suspect that its solutions are the same as
those ones listed in Theorem~\ref{T:mjor1}.
Unfortunately, this is only a conjecture which is left as an
open problem.

Finally, we point out the fact that in all of our statements we have
supposed that the underlying spaces are at least 2-dimensional. In fact,
this assumption is necessary to put as it turns out from the following
example.
Over 1-dimensional spaces, standard operator algebras are trivially
identified with the complex field $\C$.
Now, consider a bijective additive function $a:\R \to \R$ for which
$a(1)=1$ and $a(\log_2 3)\neq \log_2 3$. Such a function exists since
$1$ and $\log_2 3$ are linearly independent in the linear space $\R$
over the field of rationals.
Let $f:]0, +\infty [ \to ]0, +\infty[$ be defined by
\[
f(t)=2^{a(\log_2 t)}.
\]
Define the funtion $h:\C \to \C$ by
\[
h(z) =
\begin{cases}
0 ,                   &\text{if $z=0$;}    \\
f(|z|)\frac{z}{|z|} , &\text{if $z\neq 0$.}
\end{cases}
\]
It is easy to see that $h:\C \to \C$ is a bijective multiplicative
function, $h(2)=2$ and $h(\overline{z})=\overline{h(z)}$ $(z\in \C)$.
On the other hand, $h$ is not additive, since $h(3)\neq h(1)+h(2)$. This
function serves as a counterexample for all of our results above
after omitting the assumption on dimension.

\section{Proofs}

This section is devoted to the proofs of our results.

If $X$ is a Banach space, then the operator $P\in B(X)$ is an
idempotent if $P^2=P$. There is a partial ordering between idempotents.
If $P,Q\in B(X)$ are idempotents, then we write $P\leq Q$ if $PQ=QP=P$.
The idempotents $P,Q\in B(X)$ are said to be orthogonal if
$PQ=QP=0$.

Let $x\in X$ and $f\in X'$ be nonzero. The rank-1 operator $x\ot f$ is
defined by
\[
(x\ot f)(z)=f(z)x \qquad (z\in X).
\]
It is trivial to see that $x\ot f$ is an idempotent if and only if
$f(x)=1$. Conversely, every rank-1 idempotent can be written in this
form.

We begin with the proof of our first result.

\begin{proof}[Proof of Theorem~\ref{T:mjor1}]
First observe that $\phi(0)=0$. Indeed, if $A\in \A$ is such that
$\phi(A)=0$, then we have
\[
\phi(0)=\phi((1/2)(A0+0A))=(1/2)(\phi(A)\phi(0)+\phi(0)\phi(A))=0.
\]
We deduce from \eqref{E:mjor1} that
$\phi$ preserves the idempotents. Since $\phi^{-1}$ has the same
properties
as $\phi$, it follows that $\phi$ preserves the idempotents in both
directions.
If the idempotents $P,Q\in \A$ are orthogonal, then we have
\[
0=\phi(0)=\phi((1/2)(PQ+QP))=
(1/2)(\phi(P)\phi(Q)+\phi(Q)\phi(P)).
\]
Multiplying this equality by $\phi(Q)$ from the left and from the right
respectively, we have $\phi(Q)\phi(P)\phi(Q)=\phi(P)\phi(Q)$ and
$\phi(Q)\phi(P)\phi(Q)=\phi(Q)\phi(P)$. This implies that
\[
0=\phi(Q)\phi(P)\phi(Q)=
\phi(P)\phi(Q)=
\phi(Q)\phi(P).
\]
Therefore, $\phi$ preserves the orthogonality between idempotents
in both directions.

We assert that $\phi$ preserves the partial order $\leq$ between
the idempotents.
If $P,Q\in \A$ are idempotents and $P\leq Q$, then we obtain
\[
\phi(P)=\phi((1/2)(PQ+QP))=
(1/2)(\phi(P)\phi(Q)+\phi(Q)\phi(P)).
\]
Multiplying this equality by $\phi(Q)$ from the left and from the right
respectively, we get that $\phi(Q)\phi(P)\phi(Q)=\phi(Q)\phi(P)$ and
$\phi(Q)\phi(P)\phi(Q)=\phi(P)\phi(Q)$. This implies that
$\phi(P)=\phi(P)\phi(Q)=\phi(Q)\phi(P)$ and hence we have $\phi(P)\leq
\phi(Q)$.

It is easy to see that an idempotent $P\in \A$ is of rank $n$ if and
only
if there is a system $P_1, \ldots, P_n\in \A$ of pairwise orthogonal
nonzero
idempotents for which $P_k\leq P$ $(k=1, \ldots, n)$, but there is no
such system of $n+1$ members.
It now follows that $\phi$ preserves the rank of idempotents.

Let $P,Q\in \A$ be orthogonal finite rank idempotents. We know that
$\phi(P), \phi(Q)$ are orthogonal finite rank idempotents. As $\phi$
preserves the order, we
have $\phi(P),\phi(Q)\leq \phi(P+Q)$ implying $\phi(P)+\phi(Q)\leq
\phi(P+Q)$. Since $\phi$ preserves also the rank of idempotents,
it follows that
$\phi(P)+\phi(Q)=\phi(P+Q)$. This means that $\phi$ is orthoadditive on
the set of all finite rank idempotents in $\A$.

Let $P_1, \ldots, P_n\in \A$ be pairwise orthogonal finite rank
idempotents and $\lambda_1, \ldots, \lambda_n\in \C$. Using the
orthoadditivity of $\phi$ we have
\begin{equation}\label{E:mjor5}
\begin{gathered}
\phi(\sum_k\lambda_k P_k)=\\
\phi((1/2)((\sum_k\lambda_k P_k)(\sum_l P_l)+(\sum_l P_l)
(\sum_k\lambda_k P_k)))=\\
(1/2)(\phi(\sum_k\lambda_k P_k)\phi(\sum_l P_l)+\phi(\sum_l P_l)
\phi(\sum_k\lambda_k P_k))=\\
(1/2)(\phi(\sum_k\lambda_k P_k)\sum_l\phi(P_l)+\sum_l \phi(P_l)
\phi(\sum_k\lambda_k P_k))=\\
\sum_l (1/2)(\phi(\sum_k\lambda_k P_k)\phi(P_l)+\phi(P_l)
\phi(\sum_k\lambda_k P_k))=\\
\sum_l \phi((1/2)((\sum_k\lambda_k P_k)P_l+P_l
(\sum_k\lambda_k P_k)))=\\
\end{gathered}
\end{equation}
\begin{equation*}
\begin{gathered}
\sum_l \phi(\lambda_l P_l).
\end{gathered}
\end{equation*}

Next we prove that $\phi(-P)=-\phi(P)$ for every finite rank idempotent
$P\in \A$. Let $P$ be of rank 1.
We have
\begin{equation*}
\begin{gathered}
\phi(\lambda P)=
\phi((1/2)((\lambda P)P+P(\lambda P)))=\\
(1/2)(\phi(\lambda P)\phi(P)+\phi(P)\phi(\lambda P)).
\end{gathered}
\end{equation*}
Multiplying this equality by $\phi(P)$ from the left and from the right
respectively, we have
$\phi(P)\phi(\lambda P)\phi(P)=\phi(P)\phi(\lambda P)$ and
$\phi(P)\phi(\lambda P)\phi(P)=\phi(\lambda P)\phi(P)$.
It follows that
\[
\phi(\lambda P)=\phi(P)\phi(\lambda P)\phi(P).
\]
Since $\phi(P)$ is of rank 1, it follows from the equality above that
\begin{equation}\label{E:mjor2}
\phi(\lambda P)=\mu \phi(P)
\end{equation}
for some scalar $\mu \in \C$.
So, we obtain that $\phi(-P)=c\phi(P)$ for
some scalar $c\in \C$. Since
\[
c^2\phi(P)=(c\phi(P))^2=\phi(-P)^2=\phi((-P)^2)=\phi(P),
\]
we have $c=\pm 1$. By the injectivity of $\phi$ we get
$\phi(-P)=-\phi(P)$.
Using \eqref{E:mjor5} we deduce that
\begin{equation}\label{E:mjor17}
\phi(-P)=-\phi(P)
\end{equation}
for every finite rank idempotent $P\in \A$.

For any $A,B \in \A$ we write
\[
A\circ B=
(1/2)(AB+BA).
\]
With this notation the equation \eqref{E:mjor1} can be rewritten as
\[
\phi(A\circ B)=\phi(A)\circ \phi(B) \qquad (A,B \in \A).
\]

Let $T\in F(X)$ be arbitrary and let $P\in F(X)$ be an
idempotent.
Choose a finite rank idempotent $Q\in \A$ for which
$QT=TQ=T$ and $QP=PQ=P$.
Such a $Q$ can be constructed in the following way. Let $S\in F(X)$.
Pick a finite rank idempotent $Q^l_S$ with range containing the range
of $S$. We have $Q^l_S S=S$. Next, pick a finite dimensional subspace
$M$ of $X$ whose direct sum with the kernel $N$ of $S$ is $X$. Consider
the idempotent $Q^r_S$ with range $M$ corresponding to the direct sum
$M\oplus N=X$. We have $SQ^r_S=S$. Finally, as the partially ordered set
of all finite rank idempotents on $X$ is cofinal (see, for example,
\cite[Lemma]{Molnar3}), we can choose a finite rank idempotent
$Q$ for which $Q^l_T,Q^r_T,Q^l_P,Q^r_P \leq Q$. It is easy to check that
$Q$ has the desired properties.
Now, it requires only trivial computation to verify that
\begin{equation}\label{E:mjor4}
(2P-Q)\circ (T\circ P)=PTP.
\end{equation}
It follows that
\[
\phi(2P-Q)\circ (\phi(T)\circ \phi(P))=\phi(PTP).
\]
We prove that $\phi(2P-Q)=2\phi(P)-\phi(Q)$.
Indeed, since $Q-P$ is an idempotent which is orthogonal to $P$,
by \eqref{E:mjor5} and \eqref{E:mjor17} we can compute
\begin{equation*}
\begin{gathered}
\phi(2P-Q)=\phi(P-(Q-P))=
\phi(P)+\phi(-(Q-P))=\\
\phi(P)-\phi(Q-P)=
\phi(P)-(\phi(Q)-\phi(P))= 2\phi(P)-\phi(Q).
\end{gathered}
\end{equation*}
So, we have
\[
(2\phi(P)-\phi(Q))\circ (\phi(T)\circ \phi(P))=\phi(PTP).
\]
We assert that $\phi(Q)\phi(T)\phi(Q)=\phi(T)$ and
$\phi(Q)\phi(P)\phi(Q)=\phi(P)$. In fact, these follow from the
equalities \[
\phi(T)=
(1/2)(\phi(T)\phi(Q)+\phi(Q)\phi(T))
\]
and
\[
\phi(P)=
(1/2)(\phi(P)\phi(Q)+\phi(Q)\phi(P))
\]
after mutliplying them by $\phi(Q)$ from the left and from the right,
respectively.
Similarly as in the case of \eqref{E:mjor4}, one can now easily check
that
\[
(2\phi(P)-\phi(Q))\circ (\phi(T)\circ \phi(P))=\phi(P)\phi(T)\phi(P).
\]
Therefore, we have $\phi(PTP)=\phi(P)\phi(T)\phi(P)$.
We note that in this part of the proof we have used an idea similar to
what was followed in the proof of \cite[Lemma 1.6]{Hakeda}.

In the next section of the proof we apply some ideas from
the proof of \cite[Theorem]{Molnar}. Fix a rank-1 idempotent $P\in
\A$. By \eqref{E:mjor2}, there is a function $h_P:\C \to \C$ such that
\[
\phi(\lambda P)=h_P(\lambda) \phi(P) \qquad (\lambda \in \C).
\]
We show that $h_P$ does not depend on $P$. If $Q\in \A$ is another
rank-1 idempotent not orthogonal to $P$, then we compute
\begin{equation*}
\begin{gathered}
\phi((1/2)((\lambda P)Q+Q(\lambda P)))=
(1/2)(h_P(\lambda)\phi(P)\phi(Q)+h_P(\lambda)\phi(Q)\phi(P))=\\
h_P(\lambda)(1/2)(\phi(P)\phi(Q)+\phi(Q)\phi(P)).
\end{gathered}
\end{equation*}
We similarly have
\[
\phi((1/2)(P(\lambda Q)+(\lambda Q)P))=
h_Q(\lambda)(1/2)(\phi(P)\phi(Q)+\phi(Q)\phi(P)).
\]
Since $\phi(P)\phi(Q)+\phi(Q)\phi(P)\neq 0$ ($\phi(P)$ is not orthogonal
to $\phi(Q)$), it follows that $h_P=h_Q$.
If $Q$ is orthogonal to $P$, then we can choose a rank-1 idempotent
$R\in \A$ such that $R$ is not orthogonal to $P$ and not orthogonal to
$Q$. We have $h_P=h_R=h_Q$. Therefore, there is a function $h:\C \to \C$
such that
\begin{equation}\label{E:mjor18}
\phi(\lambda P)=h(\lambda )\phi(P)
\end{equation}
for every $\lambda \in \C$ and every rank-1 idempotent $P\in \A$.

We assert that $\phi(\lambda A)=h(\lambda )\phi(A)$ for every $A\in
F(X)$. If $A$ is a finite rank idempotent, then this follows from
\eqref{E:mjor18} and
\eqref{E:mjor5}. If $A\in F(X)$ is arbitrary, then there is a finite
rank idempotent $P$ such that $PA=AP=A$. We compute
\begin{equation*}
\begin{gathered}
\phi(\lambda A)=\phi((1/2)(A(\lambda P)+(\lambda P)A))=\\
(1/2)(\phi(A)h(\lambda)\phi(P)+h(\lambda)\phi(P)\phi(A))=
h(\lambda)\phi(A).
\end{gathered}
\end{equation*}
We next prove that $h$ is multiplicative. Let $P\in \A$ be
a nonzero finite rank idempotent. We have
\begin{equation*}
\begin{gathered}
h(\lambda \mu)\phi(P)=
\phi(\lambda \mu P)=
\phi((1/2)((\lambda P)(\mu P) +(\mu P)(\lambda P)))=\\
(1/2)(h(\lambda)\phi(P)h(\mu)\phi(P) +h(\mu)\phi(P)h(\lambda)\phi(P))=
h(\lambda )h(\mu)\phi(P)
\end{gathered}
\end{equation*}
and this shows that $h$ is multiplicative.

We prove that $h$ is additive. Let $x,y\in X$ be linearly
independent vectors, and choose linear functionals $f,g\in X'$ such
that $f(x)=1,f(y)=0$ and $g(x)=0, g(y)=1$.
Let $\lambda ,\mu \in \C$ be such that $\lambda +\mu=1$. Define
$R=(\lambda x +\mu y)\ot (f+g)$, $P =x\ot f$, $Q=y\ot g$. Clearly,
$R,P,Q$ are rank-1 idempotents and $P$ is orthogonal to $Q$.
By what we already know, we deduce
\begin{equation*}
\begin{gathered}
h(\lambda +\mu)\phi(R)=
\phi((\lambda +\mu)R)=
\phi(R(P+Q)R)=\\
\phi(R)\phi(P+ Q)\phi(R)=
\phi(R)\phi(P)\phi(R)+ \phi(R)\phi(Q)\phi(R)=\\
\phi(RPR)+ \phi(RQR)=
\phi(\lambda R)+\phi(\mu R)=
(h(\lambda )+h(\mu))\phi(R).
\end{gathered}
\end{equation*}
By the multiplicativity of $h$, we have $h(\lambda
+\mu)=h(\lambda) +h(\mu)$ whenever $\lambda +\mu \neq 0$. To see the
additivity of $h$,
it remains to prove that $h(-\lambda)=-h(\lambda)$. Since $h$ is
multiplicative, it follows that $h(-\lambda)^2=h(\lambda^2)=h(\lambda
)^2$. By the injectivity of $h$ we have the desired
equality $h(-\lambda)=-h(\lambda)$.

We now verify that $\phi$ is additive on $F(X)$.
Let $A,B\in F(X)$ be arbitrary and pick any rank-1 idempotent $P\in
F(X)$. Choose $x\in X, f\in X^*$ such that $P=x\ot f$. We
compute
\begin{equation*}
\begin{gathered}
\phi(P)\phi(A+B)\phi(P)=
\phi(P(A+B) P)=
\phi(f((A+B)x) P)=
\\
h(f((A+B)x))\phi(P)=
h(f(Ax)) \phi(P)+h(f(Bx)) \phi(P)=
\\
\phi(f(Ax) P)+\phi(f(Bx) P)=
\phi(PA P)+\phi(P B P)=
\\
\phi(P)\phi(A) \phi(P)+\phi(P) \phi(B) \phi(P)=
\phi(P)(\phi(A)+\phi(B)) \phi(P).
\end{gathered}
\end{equation*}
Since this holds true for every rank-1 idempotent $P$ on $X$, we
easily obtain that $\phi(A+B)=\phi(A)+\phi(B)$. Consequently, $\phi:
F(X) \to F(Y)$ is an additive bijection satisfying \eqref{E:mjor1}.

Since the algebra $F(X)$ (as well as every standard operator algebra) is
prime (this means
that for every $A,B \in F(X)$, the equality $AF(X) B=\{ 0\}$
implies $A=0$ or $B=0$), we can apply a result of Herstein
\cite{Herstein} to obtain that
$\phi$ is necessarily a ring isomorphism or a ring antiisomorphism
of $F(X)$.
In the isomorphic case we can apply the result in \cite{Semrl}
and obtain the desired form of $\phi$ on $F(X)$. In the finite
dimensional case we are done since in that case any standard operator
algebra coincides with $F(X)$. Observe that the argument given
in \cite{Semrl} for the finite dimensional case can be changed to give
the antiisomorphic part of our result in the finite dimensional
case.
So, let us assume that $X$ is infinite dimensional and that $\phi$ is a
ring isomorphism. By \cite{Semrl} there is a
bounded invertible either linear or conjugate-linear operator $T:X\to Y$
such that
\[
\phi(A)=TAT^{-1} \qquad (A\in F(X)).
\]
If $A\in \A$ is arbitrary, then for every finite rank idempotent $P\in
F(X)$ we have
\begin{equation*}
\begin{gathered}
(1/2)T(AP+PA)T^{-1}=
\phi((1/2)(AP+PA))=\\
(1/2)(\phi(A)\phi(P)+\phi(P)\phi(A))=
(1/2)(\phi(A)TPT^{-1}+TPT^{-1}\phi(A)).
\end{gathered}
\end{equation*}
Multiplying this equality by $T^{-1}$ from the left and by $T$ from the
right, we get
\[
AP+PA=
T^{-1}\phi(A)TP+PT^{-1}\phi(A)T.
\]
Now, multiplying this equality by $P$ from both sides, we arrive at
\[
PAP=PT^{-1}\phi(A)TP.
\]
Since $P\in \A$ was an arbitrary finite rank idempotent, it follows that
$A=T^{-1}\phi(A)T$
$(A\in \A)$. Therefore, we have $\phi(A)=TAT^{-1}$ $(A\in \A)$.

Suppose finally that $X$ is infinite dimensional and
$\phi$ is a ring antiisomorpism of $F(X)$.
Performing trivial modifications in the proofs of \cite[Theorem]{Semrl}
and \cite[Proposition 3.1]{BresarSemrl}, one can verify that
there is a bounded invertible either linear or conjugate-linear
operator $T:X' \to Y$ such that
\[
\phi(A)=TA'T^{-1} \qquad (A\in F(X)).
\]
Similarly to the isomorphic case, we can arrive at the equality
\[
P'A'+A'P'=
(AP+PA)'=
T^{-1}\phi(A)TP'+P'T^{-1}\phi(A)T.
\]
As $P'$ is an idempotent, multiplying this equality by $P'$ from both
sides, we deduce
\[
P'A'P'=P'T^{-1}\phi(A)TP'.
\]
Since, as we learn from \cite[Proposition 3.1]{BresarSemrl},
in the antiisomorphic case $X,Y$ are reflexive, it follows
that $P'$ runs
through the set of all finite rank idempotents in $B(X')$ as $P$ runs
through the set of all finite rank idempotents in $B(X)$.
So, just as in the
isomorphic case we can infer that
\[
\phi(A)=TA'T^{-1} \qquad (A\in \A).
\]
This completes the proof of the theorem.
\end{proof}

It is now easy to prove
Corollary~\ref{C:mjor3}.
We recall that the self-adjoint idempotents in $B(H)$
are called projections.

\begin{proof}[Proof of Corollary~\ref{C:mjor3}]
Clearly, Theorem~\ref{T:mjor1} can be applied.
According to that result, we have several possibilities concerning the
form  of $\phi$.
We give the proof in the case of only one such possibility. The other
cases can be handled in a quite similar way.
Suppose that $H$ is infinite dimensional.
By Theorem~\ref{T:mjor1}, we have, for example, a
bounded linear operator $T:H\to K$ such that
\[
\phi(A)=TAT^{-1} \qquad (A\in \A).
\]
Pick an arbitrary rank-1 projection $P\in \A$.
By the self-adjointness of $\phi$ we have
\[
{(T^{-1})}^*PT^*=
(TPT^{-1})^*=
TPT^{-1}.
\]
Since this holds for every rank-1 projection $P$ on $H$ we easily obtain
that the vectors ${(T^{-1})}^*x$, $Tx$ are linearly dependent for every
$x\in H$. It needs only an elementary linear algebraic argument to show
that in this case ${T^{-1}}^*$ and $T$ are necessarily linearly
dependent,
that is,  we have $T^{-1}=\lambda T^*$ for some $\lambda \in \C$.
On the other hand, it follows from \eqref{E:mjor19} that
$\phi$ sends projections to projections. This implies that the scalar
$\lambda$ above is necessarily positive. Denote $U=\sqrt \lambda
T$. We infer that $U:H\to K$ is an invertible bounded linear
operator with $U^{-1}=U^*$. This gives us that $U$ is unitary.

As for the case when $H$ is finite dimensional, we recall the
well-known fact that
if $h:\C \to \C$ is a ring automorphism of $\C$ for which
$h(\overline{\lambda})
=\overline{h(\lambda)}$ $(\lambda \in \C)$, then $h$ is either the
identity or the conjugation.
\end{proof}

The proof of Theorem~\ref{T:mjor2} will rest on the following
lemmas. Recall that an operator $A\in B(H)$ is said to be positive
if $\la Ax,x\ra \geq 0$ holds for every $x\in H$. In this case we write
$A\geq 0$.

\begin{lemma}\label{L:mjor3}
Let $H$ be a Hilbert space and $A,B\in B(H)$. Suppose that $A$ is
positive and $AB+BA=0$. Then we have $AB=BA=0$.
\end{lemma}

\begin{proof}
Since $AB=-BA$, we obtain
\[
A^2B=A(AB)=A(-BA)=(-AB)A=(BA)A=BA^2.
\]
That is, $A^2$ commutes with $B$. It is well-known that if a positive
operator $T$ commutes with an operator, then the same holds true for the
positive square root of $T$. In fact, this follows from
the fact that the square root of $T$
is the norm limit of polynomials of $T$.
Therefore, we get that $A$ commutes with $B$ which gives us that
$AB=BA=0$.
\end{proof}

In what follows let $\clrng A$ denote the closure of the range of the
operator $A\in B(H)$.

\begin{lemma}\label{L:mjor4}
Let $\A$ be a standard operator algebra on a Hilbert space. Let $A,B \in
\A$ be self-adjoint. Then we have ${\clrng A} \subset \clrng
B$ if and only if for every positive operator $C\in \A$ with $BC=0$ it
follows that $AC=0$.
\end{lemma}

\begin{proof}
All we have to do is to note that $\A$ contains all projections
of rank 1 and that the condition
${\clrng A} \subset \clrng B$ is equivalent to the condition
that $\ker B\subset \ker A$.
\end{proof}

As for the proof of our next lemma we recall the following
useful notation.
If $x,y\in H$, then $x\ot y$ stands for the operator defined by
\[
(x\ot y)(z)=\la z,y\ra x \qquad (z\in H).
\]

\begin{lemma}\label{L:mjor5}
Let $H$ be a Hilbert space.
If $A\in B(H)$ is such that $TA+AT\geq 0$ holds for every $0\leq T\in
F(H)$, then $A$ is a nonnegative scalar multiple of the identity.
\end{lemma}

\begin{proof}
First observe that $A$ is positive. Indeed, for every finite rank
projection $P$ on $H$ we have $PA+AP\geq 0$. Considering an
increasing net of finite rank projections weakly converging to the
identity, we obtain that $A+A\geq 0$ and this implies
$A\geq 0$.

If $0\neq x\in H$ is arbitrary, then we have
\[
x\ot Ax+Ax \ot x\geq 0
\]
It follows from this inequality that for any $y\in H$ we have
\[
\la y, Ax\ra \la x, y\ra +
\la y, x\ra \la Ax, y\ra \geq 0
\]
which implies that
\begin{equation}\label{E:mjor12}
\Re (\la y, Ax\ra \la x, y\ra) \geq 0.
\end{equation}
We can write $Ax=\lambda x+x^\perp$, where $\lambda \in \C$ and
$x^\perp\in H$ is a vector orthogonal to $x$. Define $y=\mu x+x^\perp$
for an arbitrary $\mu \in \C$.
It follows from \eqref{E:mjor12} that
\[
\Re (\mu \bar \lambda \| x\|^2+\| x^\perp\|^2) \bar \mu\geq 0.
\]
This implies that
\[
|\mu|^2\Re \bar \lambda \| x\|^2 +\| x^\perp \|^2 \Re \bar \mu \geq 0
\]
holds for every $\mu \in \C$. It is easy to see that we necessarily have
$\| x^\perp\|^2=0$.

The above observation yields that for every $x\in H$, the vectors $Ax$
and $x$ are
linearly dependent. As we have mentined in the proof of
Corollary~\ref{C:mjor3},
such a local linear dependence implies global linear dependence.
Therefore, it follows
that $A$ is a scalar multiple of the identity. It is clear that the
scalar in question is nonnegative.
\end{proof}

\begin{lemma}\label{L:mjor7}
Let $n\in \N$, $n>1$.
Suppose that $\psi:M_n (\C)\to M_n(\C)$ is a bijective transformation
for which
\begin{equation*}
\begin{aligned}
\psi(A^*)&=\psi(A)^*\\
\psi(AB+BA)&=\psi(A)\psi(B)+\psi(B)\psi(A)
\end{aligned}
\end{equation*}
holds for every $A,B \in M_n(\C)$.
Then $\psi$ satisfies \eqref{E:mjor19}.
\end{lemma}

\begin{proof}
First observe that $\psi$ preserves positivity in both directions.
Indeed, if $A\in M_n(\C)$ is positive, then there is a positive
$B\in M_n(\C)$ such that $2B^2=A$. We have
\[
\psi(A)=\psi(B^*B+BB^*)=
\psi(B)^*\psi(B)+\psi(B)\psi(B)^*\geq 0.
\]
Since $\psi^{-1}$ has the same properties as $\psi$, we get that
$\psi$ preserves positivity in both directions.

Let $A\in M_n(\C)$ be positive. We have
\[
\psi(A)\psi(I)+\psi(I)\psi(A)=\psi(2A)\geq 0.
\]
Since $\psi(A)$ runs through the positive elements of $M_n(\C)$, by
Lemma~\ref{L:mjor5} we infer that $\psi(I)$ is a positive scalar
multiple of the identity.
Denote $\psi(I)=\lambda I$.
Consider the transformation $\tilde \psi :M_n(\C) \to M_n(\C)$ defined
by
\[
\tilde \psi(A)=(1/\lambda)\psi(A) \qquad (A\in M_n(\C)).
\]
Since we have
\[
\psi(A)=\psi(I (A/2)+(A/2)I)=2\lambda \psi(A/2),
\]
one can easily check that $\tilde \psi$ satisfies
\[
\tilde \psi((1/2)(AB+BA))=
(1/2)(\tilde \psi(A)\tilde \psi(B)+\tilde \psi(B)\tilde \psi(A)).
\]
By Theorem~\ref{T:mjor1}, $\tilde \psi$ is additive.
It follows from the definition of $\tilde \psi$ that $\psi$ is also
additive which plainly implies the assertion of the lemma.
\end{proof}

We are now in a position to prove our final result.

\begin{proof}[Proof of Theorem~\ref{T:mjor2}]
Just as in the proof of Theorem~\ref{T:mjor1} one can prove that
$\phi(0)=0$.

We next show that $\phi$ preserves the positive elements in both
directions. This can be done quite similarly to the first part of
the proof of Lemma~\ref{L:mjor7}.

Let $A\in \A$ be positive and $B\in \A$ be arbitrary. Suppose that
$AB=BA=0$.
We have
\[
0=\phi(0)=\phi(AB+BA)=\phi(A)\phi(B)+\phi(B)\phi(A).
\]
Since $\phi(A)$ is positive, it follows from Lemma~\ref{L:mjor3} that
$\phi(A)\phi(B)=\phi(B)\phi(A)=0$.
As $\phi^{-1}$ has the same properties as $\phi$, we find that
for any two self-adjoint operators
$A,B\in \A$ one of them being positive we have $AB=0$ if and only if
$\phi(A)\phi(B)=0$.
From Lemma~\ref{L:mjor4} we deduce that for any two self-adjoint
operators $A,B \in \A$ we have $\clrng A\subset \clrng B$ if and only if
$\clrng \phi(A) \subset \clrng \phi(B)$.

Let $n$ be a positive integer. Using the spectral
theorem, one can easily
verify the following characterization of positive rank-$n$
operators.
The positive operator $A\in \A$ is of rank $n$ if and only if there
exists a system $A_1, \ldots, A_{n}\in \A$ of nonzero positive
operators such
that $\clrng A_k \subset \clrng A$ $(k=1, \ldots, n)$, $A_kA_l=0$,
$(k\neq l)$ but there is no such system of $n+1$ members.
By this characterization, $\phi$ preserves the positive rank-$n$
operators in both directions.

Let $A\in \A$ be a positive rank-$n$ operator. Then $\phi(A)$ is also
positive and is of rank $n$.
Let $B\in \A$ be any operator acting on $H_0=\clrng A$. We mean by this
that $B$ maps $H_0$ into itself and $B$ is zero on $H_0^\perp$.
One can easily verify that $\clrng (B^*B+BB^*) \subset H_0$.
Denote $C=\phi(B)$ and $K_0=\clrng \phi(A)$. It follows that
$\clrng (C^*C+CC^*)\subset K_0$. If $k\in K_0^\perp$, then we have
$k\in \ker (C^*C+CC^*)$. Since
$\la C^*Ck,k\ra +\la CC^* k,k\ra=0$, we obtain that
$\la C^*Ck,k\ra=0$ and $\la CC^*k,k\ra =0$. It follows that
$Ck=0$ and $C^*k=0$. Consequently, we get that $C(K_0^\perp)=\{ 0\}$ and
$C(K_0)\subset K_0$. Therefore, we have proved that if $B\in \A$ acts on
$\clrng A$, then $\phi(B)\in \B$ acts on $\clrng \phi(A)$.

The argument above gives us that $\phi$ sends finite rank operators to
finite rank operators. Indeed, any finite rank operator can be
considered as an
operator acting on the range of a positive finite rank operator. Since
$\phi^{-1}$ has the same properties as
$\phi$, we obtain that $\phi$ maps $F(H)$ onto $F(H)$ and, identifying
the operator algebra over $H_0$ and $K_0$ with $M_n(\C)$,
the map $\phi$ induces a bijective transformation
$\psi:M_n(\C) \to M_n(\C)$ for which
\begin{equation*}
\begin{aligned}
\psi(T^*)&=\psi(T)^*\\
\psi(TS+ST)&=\psi(T)\psi(S)+\psi(S)\psi(T)
\end{aligned}
\end{equation*}
for every $T,S \in M_n(\C)$.
Lemma~\ref{L:mjor7} tells us that $\psi$ satisfies
\eqref{E:mjor19} on $M_n(\C)$.
Since $A$ was arbitrary, it follows that $\phi$ fulfils
\eqref{E:mjor19} on $F(H)$.
Now, referring to Corollary~\ref{C:mjor3} we have the form of $\phi$ on
$F(H)$ which can be shown to be valid on the whole $\A$ in a way very
similar to
the last part of the proof of Theorem~\ref{T:mjor1}.
\end{proof}

\bibliographystyle{amsplain}

\end{document}